\documentclass[12pt]{amsart}
\usepackage{geometry} 
\usepackage{amssymb}
\usepackage{amscd}
\usepackage{amsmath}
\usepackage{amsthm}

\usepackage{mathrsfs}

\usepackage[pdfpagelabels]{hyperref}

\def\CC{\mathbb{C}}

\def\beq{\begin{equation}}
\def\eeq{\end{equation}}
\def\beqa{\begin{eqnarray}}
\def\eeqa{\end{eqnarray}}

\def\fitt{\textnormal{Fitt}}
\def\OO{\mathcal O}

\newtheorem{thm}{Theorem} [section]

\theoremstyle{definition}

\newtheorem{ex}[thm]{Example}

\newtheorem{rem}[thm]{Remark}

\begin{document}

\title{Fitting Ideals without a Presentation}

\author{Ay{\c s}e Sharland}
\address{Department of Mathematics and Applied Mathematical Sciences, University of Rhode Island, Kingston, RI 02881, US}
\email{aysharland@uri.edu}

\author{Jacob Smith}
\address{Department of Mathematics and Applied Mathematical Sciences, University of Rhode Island, Kingston, RI 02881, US}
\email{jake.smith@uri.edu}



\begin{abstract}
In this article, we investigate alternative construction of Fitting ideals of pushforward modules $f_*\mathcal{O}_{X,0}$ for finite and holomorphic map-germs from an $n$-dimensional Cohen-Macaulay space $(X,0)$ to $(\mathbb{C}^{n+1},0)$. For corank 1 map-germs, we generalize a result of D. Mond and R. Pellikaan to iteratively calculate $k$-th Fitting ideals as ideal quotients of lower ones. We also show that for a stable map-germ of any corank, the first Fitting ideal can be calculated as a quotient ideal of the Jacobian of the image and the pushforward of the ramification ideal, which is a modification of classical result of due to Piene.

\end{abstract}
\pagenumbering{arabic}

\maketitle

\section{Introduction}

Fitting Ideals, introduced by H. Fitting (\cite{fitting}), are defined for finitely generated modules as ideals generated by minors of a presentation matrix. To be more precise, let $M$ be a finitely generated module over a ring $R$, and let 
\[ R^p\xrightarrow{\varphi}  R^q \rightarrow M \rightarrow 0\]
be a presentation of $M$. Assuming the map $\varphi$ is given by a matrix $A$, the ideal $I_n(A)$ of $n\times n$-minors of $A$ is generated by all determinants of $n\times n$ submatrices of $A$. By convention, we assume $I_0(A)=A$ and $I_n(A)=0$ if $n$ is larger than the row (or column) size of $A$. The ideals $I_n(A)$ are independent of the choice of bases for the presentation (\cite[Corollary-Definition 20.4]{eisenbud}), and commute with base change (\cite[Corollary 20.5]{eisenbud}). Then, the $i$-th \emph{Fitting ideal} of $M$ is the ideal 
\[\fitt_i(M) := I_{q-i}(A).\]

Fitting ideals carry important information about modules. For example, if $M$ is a module generated by $n$ elements, then 
\[(\textnormal{Ann}(M))^n\subset \textnormal{Fitt}_0(M)\subset \textnormal{Ann}(M)\]
where $\textnormal{Ann}(M)$ denotes the annihilator of $M$, in particular, their radicals are equal $\sqrt{\textnormal{Fitt}_0(M)}=\sqrt{\textnormal{Ann}(M)}$ (\cite[Proposition 10.7]{eisenbud}). For certain modules -- as we will soon see examples below, $ \textnormal{Fitt}_0(M)$ coincides with $\textnormal{Ann}(M)$. For more applications and examples, we refer readers to \cite{eisenbud}, \cite{bruns-herzog}. 

In this article, we will focus on Fitting ideals in the theory of finite and smooth differentiable map-germs between complex spaces where they appear as ideals defining multiple point spaces. Our motivation for alternative constructions lies in the fact that the size of matrices or the number of variables involved sometimes make calculations quite difficult.

In the next section we recall some terminology and notations for our context, and then cover Fitting ideals in ``order", starting with known results from literature regarding the zero-th and first Fitting ideals and then showing our main results, Theorem \ref{thm1} and Theorem \ref{thm2}. 
Throughout the paper, we will demonstrate the results with example codes in {\sc{Singular}} (\cite{singular}) and \textit{Macaulay2} (\cite{M2}). We also include an appendix to show a calculation following the traditional approach for finding the Fitting ideals for our main example throughout.

\section{Terminology and Notation}

Here we give a brief summary of terminology and notations. For details, one can refer to \cite{wall}, \cite{martinet} for the theory of map-germs and \cite{marar-mond}, \cite{mond-pellikaan} for the main results for multiple point spaces, or more recent book \cite{mond-nuno} by D. Mond and J. J. Nu{\~n}o-Ballesteros for both aspects.

One of the important studies in the theory of differentiable smooth map-germs $(\CC^n,0)$ to $(\CC^p,0)$ is the classification of map-germs. Therefore they are studied with respect to equivalence relations. Among those relations, $\mathcal{A}$-\textit{equivalence} carries the richest structure for map-germs.   

We say that two map-germs $f,g\colon (\mathbb{C}^n,0) \rightarrow (\mathbb{C}^{p},0)$ are $\mathcal{A}$-equivalent if there exist local diffeomorphisms $\varphi$ of $(\mathbb{C}^n,0)$ and $\psi$ of $(\mathbb{C}^{p},0)$ such that $g=\psi\circ f \circ \varphi^{-1}$.

Assume that $f$ has components $f=(f_1,f_2,\ldots, f_{n+1})$ with $f_i\colon (\mathbb{C}^n,0)\rightarrow (\mathbb{C},0)$. The $k$-\textit{jet} $j^k(f)$ is the polynomial map given by $j^k(f)=(j^k(f_1),\ldots, j^k(f_{n+1}))$ where $j^k(f_i)$ is the Taylor polynomial of $f_i$ of degree $k$ (without a constant term).
Then, a map-germ $f$ is $k$-\textit{determined} if every map-germ $g\colon (\mathbb{C}^n,0) \rightarrow (\mathbb{C}^{p},0)$ with the same $k$-jet is $\mathcal{A}$-equivalent $f$. We say that $f$ is finitely $\mathcal{A}$-determined (or $\mathcal{A}$-finite) if it is $k$-determined for some $k\leq \infty$. 

A map-germ is $\mathcal{A}$-\textit{stable} if any of its unfoldings is $\mathcal{A}$-equivalent to the trivial unfolding $f\times 1$. By fundamental results of Mather (\cite{matherII}), finite determinacy is equivalent to the finite dimensionality of the normal space 
\begin{equation}\label{normalsp}N\mathcal{A}_ef:=\frac{f^*(\Theta_{\mathbb{C}^p,0})}{tf(\Theta_{\mathbb{C}^n,0})+f^{-1}(\Theta_{\mathbb{C}^p,0})},\end{equation} and thus (if $f$ is not stable) to $0\in\mathbb{C}^p$ being an isolated point of instability of $f$. We set  $\mathcal{A}_e\textnormal{-codim}(f):=\textnormal{dim}_{\mathbb{C}}N\mathcal{A}_ef$. 

We will denote the ring of holomorphic functions on a germ $(X,0)$ by $\OO_{X,0}$ and its maximal ideal by $\mathfrak{m}_{X,0}$. When needed, we will use variables $\mathbf{x}=(x_1,\ldots,x_n)$ for local coordinates in the source and $\mathbf{Y}=(Y_1,\ldots,Y_{n+1})$ in the target space for a map-germ.

The corank of a map-germ $f\colon (\mathbb{C}^n,0) \rightarrow (\mathbb{C}^{p},0)$ with $n\leq p$ is defined to be 
\[ \textnormal{corank }f=n-\textnormal{rank } \textnormal{d}f(0),\]
and its \textit{ramification} ideal $R_f$ is defined to be 
\[ R_f:=\textnormal{Fitt}_0(\Omega_f)\] 
i.e. the ideal generated by all $n\times n$-minors of $\textnormal{d}f$. On the other hand, if $X\subset \CC^n$ is a $d$-dimensional variety defined an ideal  $I=(f_1,\ldots,f_s)$, then the Jacobian ideal of $X$ is the ideal of $(n-d)\times(n-d)$-minors of the matrix $\left(\frac{\partial f_i}{\partial x_j}\right)_{1\leq i \leq s \atop 1\leq j \leq n}$.

For a finite map $f\colon X \rightarrow Y$, \textit{$k$-th multiple point space} is defined to be the set $M_k(f)$ of points in $Y$ whose preimage consists of $k$ or more points, counting multiplicity. If $X$ and $Y$ are analytic spaces, then $\mathcal{O}_{X}$ is a finite module over $\mathcal{O}_Y$ via $f$, and 
there exists a presentation of $f_*\mathcal{O}_{X}$ over $\mathcal{O}_Y$ from which one can calculate the Fitting ideals $\fitt_i(f_*\OO_{X,0})$. Furthermore, the $k$-th Fitting ideal gives $M_{k+1}(f)$ an analytic structure; that is, $M_{k+1}(f)$ is the zero set $V(\textnormal{Fitt}_k(f_*\mathcal{O}_{X}))$ (\cite[Proposition 1.5]{mond-pellikaan}). 
For example, the $0$-th Fitting Ideal defines the image of $f$, $V(\textnormal{Fitt}_1( f_*\mathcal{O}_{X})) $ consists of double points and ramification points, $V(\textnormal{Fitt}_2( f_*\mathcal{O}_{X}))$ is the set of triple points and so on.

In the case that $(X,0)$ is a germ of an irreducible Cohen-Macaulay space of dimension $n$ and $Y\cong \mathbb{C}^{n+1}$, the presentation of $f_*\mathcal{O}_{X}$ is a short exact sequence
\[0\rightarrow \mathcal{O}_{\mathbb{C}^{n+1},0}^{r+1}\xrightarrow{\lambda}  \mathcal{O}_{\mathbb{C}^{n+1},0}^{r+1} \xrightarrow{G} f_*\mathcal{O}_{X}\rightarrow 0\]
where $\lambda$ is a symmetric matrix and $G$ is induced by the generators $g_0,g_1,\ldots, g_r$ of the $Q$-\textit{algebra} of $f$, $Q(f):=\mathcal{O}_{\CC^n,0}/f^*{\mathfrak{m}_{\CC^{n+1},0}}$.  An algorithm for constructing such presentation, namely the matrix of $\lambda$, is given by Mond and Pellikaan in \cite[Section 2.2]{mond-pellikaan}. We will visit their algorithm in Appendix \ref{app1} for our main example in this paper.

It is also advantageous to consider $k$-th multiple point spaces $D^k(f)$ over the source:  
\[D^k(f)=\textnormal{closure}\{(\mathbf{x}_1,\ldots,\mathbf{x}_k)\in X^k| f(\mathbf{x}_1)=\cdots=f(\mathbf{x}_k), \mathbf{x}_i\neq \mathbf{x}_j \text{ if } i\neq j\}.\]
There are no explicit formula for ideals defining $D^k(f)$ in the general setting, accept for $k=2$. The ideal
\begin{equation}\label{d2} \mathcal{I}_2:=(f\times f)^*I(\Delta_p)+\textnormal{Fitt}_0(
I(\Delta_n)/(f\times f)^*I(\Delta_p))
\end{equation}
where $I(\Delta_n)$ and $I(\Delta_p)$ are the ideal sheaves defining the diagonals 
$\Delta_n$ in $\mathbb{C}^n\times \mathbb{C}^n$ and $\Delta_p$ in $\mathbb{C}^p\times 
\mathbb{C}^p$, gives $D^2(f)$ an analytic structure. On the other hand, if the corank of $f$ is equal to 1, an explicit list of generators for the ideal defining $D^k(f)$ in $(\CC^n)^k$ is available for all $k\geq 2$ and can be found in \cite{mond87} and \cite{marar-mond}.  

We also consider the natural projections $\pi^{k+1}_{k,i}\colon D^{k+1}(f)\rightarrow D^k(f)$ induced by the projection $(\mathbb{C}^n)^{k+1}\rightarrow (\mathbb{C}^n)^{k}$ which forgets the $i$-th factor. The images of  $\pi^{k+1}_{k,i}$ for different $i$ and the same $k$ are equal to each other. We take $i=k+1$ and simply write $\pi^{k+1}_k$ instead of $\pi^{k+1}_{k,i}$, and set $\pi^1_0:=f$. Moreover, let $D^k_j(f)=\pi^{j+1}_{j}\circ \cdots \circ \pi^{k}_{k-1}(D^k(f))$ and  $\epsilon^k=f\circ \pi^2_1\circ \cdots \pi^k_{k-1}$. Then we have $\epsilon^k(D^k(f))=M_k(f)$ and 
$D^k_1(f)=f^{-1}(M_k(f))$ as sets.

\section{The Zero-th Fitting Ideal} 

Let $(X,0)$ be an $n$-dimensional Cohen-Macaulay germ of an irreducible variety. Let $f\colon (X,0)\rightarrow (\mathbb{C}^{n+1},0)$ be a generically one-to-one 
map-germ whose image $V$ is defined by an ideal $I\subset \mathcal{O}_{\mathbb{C}^{n+1},0}$. Then, for example by \cite[Proposition 3.1]{mond-pellikaan}, $\textnormal{Fitt}_0(f_*\mathcal{O}_{X,0})= I$ and is a principal ideal. Therefore, it can be calculated by finding the minimal relation between the components of $f$.


\begin{ex}\label{ex-elim} Let us consider the map-germ 
\begin{eqnarray} \label{ex0-map} f\colon (\CC^2,0)&\rightarrow &(\CC^3,0)\\
(x,y)&\mapsto &(x,y^5-xy,y^6+xy^2)\nonumber \end{eqnarray}
which is finite and weighted homogeneous with weights $(4,1)$ and degrees $(4,5,6)$. It has $\mathcal{A}_e$-codimension 20 which can be easily calculated by Mond's formula given in \cite[Theorem 1]{mond-novan} (for an alternative method, based on (\ref{njg}), see \cite{altintas-ex3}, also \cite[Algorithm 8.1]{mond-nuno}).
Let us denote the coordinates on $(\CC^3,0)$ by $(X,Y,Z)$ for simplicity. The following {\sc{Singular}} code demonstrates the elimination method for finding the ideal $I=(h)$ defining the image of $f$. 

\begin{verbatim}
ring t=0,(X,Y,Z),(wp(4,5,6));
ring s=0,(x,y),(wp(4,1));
def st=s+t;
setring st;
ideal I1=X-x,Y-(y5-xy),Z-(y6+xy2);
ideal h=eliminate(I1,xy);
h;
->
h[1]=16X5Y2+Y6-16X6Z+11XY4Z+28X2Y2Z2+8X3Z3-Z5
\end{verbatim}
Alternatively, \textsf{preimage} command can be used for calculating the image of a map-germ.

\begin{verbatim}
ring t=0,(X,Y,Z),(wp(4,5,6));
ring s=0,(x,y),(wp(4,1));
ideal p=0;
map f=t,x,y5-xy,y6+xy2;
setring t;
ideal h=preimage(s,f,p);
\end{verbatim}

For weighted homogeneous examples, such as our current one, \textit{Macaulay2} provides another neat option.

\begin{verbatim}
S=QQ[x,y,Degrees=>{4,1}]
T=QQ[X,Y,Z,Degrees=>{4,5,6}]
f=map(S,T,{x,y^5-x*y,y^6+x*y^2})
h=ker f
              5 2    6      6         4       2 2 2     3 3    5
o4 = ideal(16X Y  + Y  - 16X Z + 11X*Y Z + 28X Y Z  + 8X Z  - Z )

o4 : Ideal of T
\end{verbatim}
\end{ex}

\section{The First Fitting Ideal}

In this section, we note a couple of known results from the literature for constructing the first Fitting Ideal without a presentation. And then show our first main result in Section \ref{sect-piene} motivated by Piene's main result in \cite{piene79}.
First we recall the definition of the common ingredient for this section: the \textit{conductor} ideal and its connection to the first Fitting Ideal.

Let $f\colon (X,0) \rightarrow (Y,0) $ be a finite map-germ of analytic spaces. 
The \textit{conductor ideal} $\mathcal{C}$ of $X$ over the image $V$ of $f$ is defined to be \(\mathcal{C}:=\textnormal{Hom}_{\mathcal{O}_{V,0}}(\mathcal{O}_{X,0} ,\mathcal{O}_{V,0} ),\)
that is,
\[\mathcal{C}=\{ g  \in \mathcal{O}_{V,0} \mid g\cdot \mathcal{O}_{X,0} \subset \mathcal{O}_{V,0}\}\]
(see, for example, \cite{jong-straten-norm}, \cite{mond-pellikaan}, \cite{mond-nuno},\cite{sing-book}).  Since $\OO_{X,0}$ has a $\OO_{Y,0}$-module structure induced by $f$, $\mathcal{C}$ can  be seen as an ideal of both $\OO_{X,0}$ and $\OO_{V,0}$. It is, in fact, the largest of such ideals. 

By \cite[Proposition 3.5]{mond-pellikaan}, if $(X,0)$ is a Cohen-Macaulay space of dimension $n$ and $f\colon (X,0)\rightarrow (\CC^{n+1},0)$ is finite and generically one-to-one,
\begin{equation}\label{fitt1-cond} \fitt_1(f_*{\mathcal{O}_{X,0}})\OO_V = \mathcal{C}, \end{equation}
which yields $(f^*)^{-1}f^*\mathcal{C}=\fitt_1(f_*{\mathcal{O}_{X,0}})$ 
(also see \cite[Proposition 3.4]{klu-curv}, \cite[Lemma 11.1]{mond-nuno}). In what follow, we will use the conductor ideal to determine $\fitt_1(f_*{\mathcal{O}_{X,0}})$.

 

\subsection{As the radical of the Jacobian of the image.} Let us assume that  $f\colon (\CC^n,0)\rightarrow (\mathbb{C}^{n+1},0)$ a generically one-to-one map-germ whose image is $V=V(h)$ for some $h\in \OO_{\CC^{n+1},0}$. Then, the singular locus of the image of $f$ is defined by the Jacobian ideal
\[J_h:=\left(\frac{\partial h}{\partial Y_1},\ldots, \frac{\partial h}{\partial Y_{n+1}} \right)\]
In general, the Jacobian ideal is not reduced. On the other hand, the conductor ideal $\mathcal{C}$, as an ideal of $\OO_{\CC^{n+1},0}$, is the ideal of functions vanishing on the singular locus of the image is reduced, provided that $f$ is $\mathcal{A}$-finite (\cite[Proposition 8.4(3)]{mond-nuno}), whence $\mathcal{C}=\sqrt{J(h)}$. 

\begin{ex}\label{ex-rad} Let $f$ be the map-germ we studied in Example \ref{ex-elim}. We can calculate $\fitt_1(f_*\OO_{\CC^n,0})$ by the following code and obtain the desired result (cf. Example \ref{ex0}). 

\begin{verbatim}
LIB "ring.lib";
ring t=0,(X,Y,Z),(wp(4,5,6));
ring s=0,(x,y),(wp(4,1));
ideal p=0;
map f=t,x,y5-xy,y6+xy2;
setring t;
ideal h=preimage(s,f,p);
ideal jh=jacob(h);
ideal fit1=radical(jh);
fit1;
->
fit1[1]=Y5+7XY3Z+8X2YZ2
fit1[2]=2XY4+10X2Y2Z+4X3Z2-Z4
fit1[3]=4X2Y3+12X3YZ+YZ3
fit1[4]=8X3Y2+8X4Z-Y2Z2-2XZ3
fit1[5]=16X4Y+Y3Z+4XYZ2
fit1[6]=16X5-Y4-6XY2Z-4X2Z2
\end{verbatim}
\end{ex}

\subsection{As the non-normal locus} 
Let $A$ be a Noetherian ring and let $N(A)$ be the \textit{non-normal} locus of $A$, that is, the set of points $p\in \textnormal{Spec}(A)$ such that $A_p$ is not normal. By, for example, \cite[Lemma 3.6.3]{sing-book}, $N(A)$ is the zero set of the conductor of $A$ in its integral closure $\bar{A}$. 

In our setting, a finite and generically one-to-one map-germ $f\colon (\CC^n,0) \rightarrow (\CC^{n+1},0)$, with image $V=V(h)$, is naturally a normalization of $V$. Therefore we can recover $\mathcal{C}$ as using the procedure \textsf{normalConductor} from \textsf{normal.lib} in {\sc{Singular}}.

\begin{ex} We consider the map-germ in Example \ref{ex-elim} to compare with the results of Example \ref{ex-rad} and \ref{ex0}.
\end{ex}
 
\begin{verbatim}
LIB "normal.lib";
ring t=0,(X,Y,Z),(wp(4,5,6));
ring s=0,(x,y),(wp(4,1));
ideal p=0;
map f=t,x,y5-xy,y6+xy2;
setring t;
ideal h=preimage(s,f,p);
ideal fit1=normalConductor(h);
fit1;
->
fit1[1]=16X5-Y4-6XY2Z-4X2Z2
fit1[2]=16X4Y+Y3Z+4XYZ2
fit1[3]=8X3Y2+8X4Z-Y2Z2-2XZ3
fit1[4]=4X2Y3+12X3YZ+YZ3
fit1[5]=2XY4+10X2Y2Z+4X3Z2-Z4
fit1[6]=Y5+7XY3Z+8X2YZ2

\end{verbatim}

\begin{rem} \textit{Grauert-Remmert criterion} for normality (\cite{grauert-remmert-coherent}) states that $A$ is normal if and only if $A=\textnormal{Hom}_A(J,J)$ for a radical ideal $J$ with $N(A)\subset V(J)$ and containing a non-zero divisor of $A$. Following this criteria, in \cite[Corollary 3.6.11]{sing-book}, Greuel and Pfister show that $N(A)$ is defined by $\textnormal{Ann}_A(\textnormal{Hom}_A(J,J)/A)$ and that 
\begin{equation}\label{gr} \textnormal{Ann}_A(\textnormal{Hom}_A(J,J)/A) =  ((a) : (aJ:J))\end{equation} for any non-zerodivisor $a\in J$. They also note that $J$ is can be chosen to be the ideal of functions vanishing on the singular locus of $A$. In our examples, we observe that working with $J_h$ rather than its radical still returns the correct ideal for $\mathcal{C}$ as we show below (cf. \cite[Example 3.6.13]{sing-book}). Although we have not found an argument to drop the requirement that $J$ has to be radical in (\ref{gr}).

\begin{verbatim}
ring t=0,(X,Y,Z),(wp(4,5,6));
ring s=0,(x,y),(wp(4,1));
ideal p=0;
map f=t,x,y5-xy,y6+xy2;
setring t;
ideal h=preimage(s,f,p);
ideal jh=jacob(h);
ideal a=jh[1];
qring q=std(h);
ideal jh=imap(t,jh);
ideal a=imap(t,a);
ideal fit1=quotient(a,quotient(a*jh,jh));
fit1;
->
fit1[1]=16X5-Y4-6XY2Z-4X2Z2
fit1[2]=16X4Y+Y3Z+4XYZ2
fit1[3]=8X3Y2+8X4Z-Y2Z2-2XZ3
fit1[4]=4X2Y3+12X3YZ+YZ3
fit1[5]=2XY4+10X2Y2Z+4X3Z2-Z4
fit1[6]=Y5+7XY3Z+8X2YZ2
\end{verbatim}

\end{rem}

\subsection{Carrying Piene's result over the target}\label{sect-piene}

In \cite{piene79}, Piene showed an important relation between the conductor and the Jacobian ideal for normalizations. To be more precise, let $f\colon X\rightarrow Z$ be a normalization (or a resolution of singularities), and $Z$ a local complete intersection with the Jacobian ideal $J_Z$. Then 
\begin{equation} \label{eq-piene} f^*J_Z=R_f\cdot f^*\mathcal{C}. \end{equation}
Consequently, we have 
\begin{equation} \label{eq-piene2} f^*\mathcal{C}=(f^*J_Z:R_f). \end{equation}
An alternative proof for (\ref{eq-piene2}) for map-germs $f\colon (\CC^{n},0) \rightarrow (\CC^{n+1},0)$, whose ramification locus has codimension 2, is given by Bruce and Marar in \cite[Proposition 2.4]{bruce-marar}.

Our result below shows a relation similar to (\ref{eq-piene2}) over the target ring for stable map-germs.

\begin{thm}[\cite{jake-thesis}]\label{thm1} Let $F\colon (\CC^{N},0) \rightarrow (\CC^{N+1},0)$ be a weighted-homogeneous stable map-germ. Then 
\begin{equation}\label{eq-thm1} (J_H:(F^*)^{-1}(R_F))=\fitt_1(F_*\mathcal{O}_{\CC^N,0}). \end{equation}
\end{thm}

\begin{proof} Our proof is mainly based on a couple of important results due to Mond (\cite{mond89}): the isomorphism 
\begin{equation} \label{njg} \frac{\theta(f)}{T\mathcal{A}_e f}\cong \frac{J_h\mathcal{O}_{\CC^{n},0}}{J_h\mathcal{O}_{V,0}} \end{equation} for any finitely $\mathcal{A}$-determined map-germ $f$, and the short exact sequence of $\OO_{\CC^{n+1},0}$-modules\begin{equation}\label{mond-ses} 0\rightarrow \frac{J_h\mathcal{O}_{\CC^{n},0}}{J_h\mathcal{O}_{V_0}} \rightarrow \frac{\mathcal{C}}{J_h\mathcal{O}_{V,0}} \rightarrow \frac{\mathcal{C}}{J_h\mathcal{O}_{\CC^{n},0}} \rightarrow  0\end{equation}  which he used to relate the topology of a stable perturbation of $f$ to its $\mathcal{A}_e$-codimension.

Firstly, if $F$ is stable then it is also infinitesimally stable, that is, $N\mathcal{A}_eF$ is trivial. Then, by (\ref{njg}), the short exact sequence (\ref{mond-ses}), set up for $F$, reduces down to
\begin{equation}\label{jg-isom} \frac{\mathcal{C}}{J_H\mathcal{O}_{V,0}} \cong \frac{\mathcal{C}}{J_H\mathcal{O}_{\CC^{N},0}} \end{equation}
where $V$ is the image of $F$ defined to be $V=V(H)$ and $J_H$ is the Jacobian ideal. 
By Piene's result (also see \cite[Proposition 2.3]{mond89}), there is an isomorphism 
\[\frac{\mathcal{O}_{\mathbb{C}^N,0}}{R_F}\cong \frac{\mathcal{C}}{J_H\mathcal{O}_{\CC^{N},0}}.\]
It follows from (\ref{jg-isom}) that $\mathcal{C}/J_H\mathcal{O}_{V,0}\cong \mathcal{O}_{\mathbb{C}^{N},0}/R_F$ over $\mathcal{O}_{\mathbb{C}^{N+1},0}$. Therefore, $(F^{*})^{-1}R_F$ annihilates $\mathcal{C}/J_H\mathcal{O}_{V,0}$, and that $\mathcal{C}\subseteq (J_H \mathcal{O}_{V,0}: (F^{*})^{-1}R_F)$. Since $F$ is weighted homogeneous, $J_H \mathcal{O}_{V,0}=J_H$. So,
\begin{equation}\label{cinc} \mathcal{C}\subseteq (J_H: (F^{*})^{-1}R_F) .\end{equation}
In order to show the other inclusion, let us consider an element $b\in (J_H: (F^{*})^{-1}R_F)$. Then, by definition, $b\cdot \gamma \in J_H$ for any $\gamma\in (F^{*})^{-1}R_F$. Applying the pullback $F^*$,
\[F^*(b\cdot \gamma)=F^*b\cdot F^*\gamma \in F^*J_H.\]
Again, by Piene's result, $F^*J_H=F^*\mathcal{C}\cdot R_F$. Also note that $F^*\gamma\in F*(F^{*})^{-1}R_F$ and that $F^*(F^{*})^{-1}R_F\subset R_F$ since $R_F$ is reduced and Cohen-Macaulay in the case of stable map-germs. Therefore, we must have $F^*b\in F^*\mathcal{C}$, in return, $b\in \mathcal{C}$. This completes the proof.

\end{proof}

If $F$ is a stable unfolding of a finite map-germ $f\colon (\CC^n,0)\rightarrow (\CC^{n+1},0)$ with $d$-parameters, that is, $F\colon (\CC^n\times \CC^d,0)\rightarrow (\CC^{n+1}\times \CC^d,0)$ given by $F(\mathbf{x},\mathbf{u})=(F_\mathbf{u}(\mathbf{x}),\mathbf{u})$ with $F_0(\mathbf{x})=f(\mathbf{x})$, then we have
\[\fitt_1(f_*\OO_{\CC^n,0}) = \fitt_1(F_*\OO_{\CC^n\times\CC^d,0})\otimes_{\OO_{\CC^{n+1}\times\CC^d,0}}\OO_{\CC^{n+1}\times\CC^d,0}/\mathfrak{m}_{\CC^d,0}\]
since Fitting ideals commute with base change. Although, we should note that calculating $\fitt_1(F_*\OO_{\CC^n\times\CC^d,0})$ this way may not be the quickest way of calculations since it will naturally involve more variables. Regardless, our result might open doors to other unsolved problems in the theory (\cite{mond-open}, \cite{mond-nuno}).

\begin{ex}  Let us consider the the map-germ $f$ in (\ref{ex0-map}). A stable unfolding $F$ is given by
\begin{eqnarray*} F\colon (\CC^2\times \CC^6,0) & \rightarrow & (\CC^3\times \CC^6,0)\\
(x,y,a,b,c,u,v,w)&\mapsto& (x,y^5+x^2y+ay+by^2+cy^3,y^6+xy^4+uy+vy^2+wy^3,a,b,c,u,v,w)
\end{eqnarray*}
based on Mather's theory, \cite[Lemma 5.9, Theorem 5.10]{matherIV}.
\begin{verbatim}
LIB "ring.lib"; 
ring t=0,(X,Y,Z,A,B,C,U,V,W),(wp(4,5,6,2,3,4,2,3,5));
ring s=0,(x,y,a,b,c,u,v,w),(wp(4,1,2,3,4,2,3,5));
ideal p=0;
map f=t,x,y5-xy+ay3+by2+cy,y6+xy2+uy4+vy3+wy,a,b,c,u,v,w;
ideal g=x,y5-xy+ay3+by2+cy,y6+xy2+uy4+vy3+wy,a,b,c,u,v,w;
matrix jF=jacob(g);
ideal RF=std(minor(jF,8));
setring t;
ideal h=preimage(s,f,p);
ideal jh=jacob(h);
ideal FRF=preimage(s,f,RF);
ideal fit1F=quotient(jh,FRF);
\end{verbatim}

Then using
\[\fitt_1(f_*\OO_{\CC^2,0}) = \fitt_1(F_*\OO_{\CC^2\times\CC^6,0})\otimes \OO_{\CC^3\times\CC^6,0}/\mathfrak{m}_{\CC^6,0}\]
and continuing the code above, we find the same result for the first Fitting Ideal as we have in Examples \ref{ex-rad} and \ref{ex0}.
\begin{verbatim}
ring t0=0,(X,Y,Z),(wp(4,5,6));
ideal fit1=std(imap(t,fit1F));
\end{verbatim}

\end{ex}

\begin{rem} 
Once we have $h$ and the first Fitting ideal, we can recover the presentation matrix $\lambda$ following the arguments (and proofs) of Propositions 3.4 and 3.14 in \cite{mond-pellikaan}: Let $\lambda$ be $(r+1)\times(r+1)$-matrix associated to a finite and generically one-to-one map-germ $f\colon (X,0)\rightarrow (\CC^{n+1},0)$ and $\lambda^1$ the matrix obtained from $\lambda$ by deleting the first row. Then $\fitt_1(f_*\OO_{X,0})=I_r(\lambda^1)$, and the sequence
\[\OO_{\CC^{n+1},0}^{r}\xrightarrow{(\lambda^1)^t} \OO_{\CC^{n+1},0}^{r+1} \rightarrow \fitt_1(f_*\OO_{X,0}) \rightarrow 0\]
is exact. Moreover, by considering the Laplace expansion for $h=\textnormal{det}(\lambda)$ in terms of the first row of $\lambda$,
\[\OO_{V,0}^{r+1}\xrightarrow{\lambda^t} \OO_{V,0}^{r+1} \rightarrow \fitt_1(f_*\OO_{X,0})\OO_{V,0} \rightarrow 0\]
is also exact. Therefore, we can find $\lambda$ matrix by calculating the first syzygy of $\fitt_1(f_*\OO_{X,0})\OO_{V,0}=\mathcal{C} $. The following code shows an example in Singular.\end{rem}

\begin{ex}
\begin{verbatim}
LIB "ring.lib";
ring t=0,(X,Y,Z),(wp(4,5,6));
ring s=0,(x,y),(wp(4,1));
ideal p=0;
map f=t,x,y5-xy,y6+xy2;
setring t;
ideal h=preimage(s,f,p);
ideal jh=jacob(h);
ideal fit1=radical(jh);
qring q=std(h);
ideal fit1=imap(t,fit1);
list resfit1=mres(fit1,0);
def lambda0=resfit1[2];  
print(lambda0);  
->
Z,  0,  0,  2XY, Y2,  
Y,  Z,  0,  -2X2,-3XY,
-2X,Y,  Z,  0,   2X2, 
0,  -2X,Y,  Z,   0,   
0,  0,  -2X,Y,   -Z  
\end{verbatim}
To see that the output is equivalent to the matrix $\lambda$ in (\ref{ex0-lambda}), we need the following matrix operations.
\begin{verbatim}
def pr1=multcol(lambda0,1,-1);
def pr2=multcol(pr1,2,-1);
def pr3=multcol(pr2,3,-1);
def pr4=multcol(pr3,4,-1);
map i=q,X,-Y,Z;
def pr5=i(pr4);
def lambda=addcol(pr5,1,-X,5);
\end{verbatim}

\end{ex}
\section{Higher Order Fitting Ideals}

Our result in this section is motivated by Mond and Pellikaan's construction of the second Fitting ideal as a quotient in \cite[Proposition 4.1]{mond-pellikaan}, namely, that
\begin{equation}\label{fit2} (\textnormal{Fitt}_1^2(f_*\mathcal{O}_{X,0}) : \textnormal{Fitt}_0(f_*\mathcal{O}_{X,0}))=\textnormal{Fitt}_2(f_*\mathcal{O}_{X,0})
\end{equation}
where $\textnormal{Fitt}_1^2(f_*\mathcal{O}_{X,0})=\textnormal{Fitt}_1(f_*\mathcal{O}_{X,0})\cdot \textnormal{Fitt}_1(f_*\mathcal{O}_{X,0})$.

The following iteration of (\ref{fit2}) was proved for stable map-germs of corank 1 in \cite{jake-thesis}. Here, we extend the result to finitely $\mathcal{A}$-determined map-germs of corank 1.
\begin{thm} \label{thm2}
	Let $X$ be an $n$-dimensional Cohen-Macaulay space and $f\colon (X,0) \to (\CC^{n+1},0)$ an $\mathcal{A}$-finite map-germ of corank 1. Also assume that $f_*\mathcal{O}_{X,0}$ has a presentation over $\mathcal{O}_{\mathbb{C}^{n+1},0}$ given by an $(r+1)\times (r+1)$-matrix $\lambda$.
	Then for each $0\leq k\leq r-1$, we have
	\begin{equation}\label{thm-fitk}
		(\fitt_{k+1}^2(f_*\mathcal{O}_{X,0}) : \fitt_{k}(f_*\mathcal{O}_{X,0})) 
		= \fitt_{k+2}(f_*\mathcal{O}_{X,0}).
	\end{equation}
\end{thm}

\begin{proof} Firstly, notice that our claim follows from the definition of Fitting ideals for $k=r-1$: we have $\fitt_{k+2}(f_*\mathcal{O}_{X,0})=1$ and $\fitt_{k+1}(f_*\mathcal{O}_{X,0})$ generated the entries of $\lambda$. It is clear that $\fitt_{k}(f_*\mathcal{O}_{X,0})$, ideal generated by $2\times 2$-minors of $\lambda$, is contained in  $\fitt_{k+1}^2(f_*\mathcal{O}_{X,0})$. 

For other cases, we will prove the result by an induction on $k$. The initial case $k=0$ is (\ref{fit2}). Let us assume that the claim (\ref{thm-fitk}) holds for some $k\leq r-2$. We want to show that
\[
		(\fitt_{k+2}^2(f_*\mathcal{O}_{X,0}) : \fitt_{k+1}(f_*\mathcal{O}_{X,0})) 
		= \fitt_{k+3}(f_*\mathcal{O}_{X,0}).
\]
The inclusion \[\fitt_{k+3}(f_*\mathcal{O}_{X,0})
		\subseteq (\fitt_{k+2}^2(f_*\mathcal{O}_{X,0}) : \fitt_{k+1}(f_*\mathcal{O}_{X,0}))
\]
follows from a special case of \cite[Lemma 10.10]{bruns-vetter} where it is shown that the product of any $u\times u$ and $v\times v$ minors of a matrix belongs, up to a scale, to the $\mathbb{Z}$-module generated by $\delta_{u+1}\delta_{v-1}$ where $\delta_{u+1}$ (and $\delta_{v-1}$ resp.) vary over all $(u+1)\times (u+1)$-minors (and $(v-1)\times (v-1)$-minors).

For the other inclusion, let us take $\beta\in (\fitt_{k+2}^2(f_*\mathcal{O}_{X,0}) : \fitt_{k+1}(f_*\mathcal{O}_{X,0}))$. In order to show that $\beta$ is an element of $\fitt_{k+3}(f_*\mathcal{O}_{X,0})$, we will relate the Fitting ideals of $f$ to those of $\pi^2_1\colon D^2(f)\rightarrow (X,0)$. 

By the definition of quotient of ideals, we have
\[\beta \cdot \fitt_{k+1}(f_*\mathcal{O}_{X,0}) \subseteq \fitt_{k+2}^2(f_*\mathcal{O}_{X,0}). \]
If we apply the pullback map $f^*$, we get
\begin{equation} \label{fs-inc} f^*(\beta \cdot \fitt_{k+1}(f_*\mathcal{O}_{X,0})) \subseteq f^*\fitt_{k+2}^2(f_*\mathcal{O}_{X,0}). \end{equation}
Since $f^*$ is a ring homomorphism, (\ref{fs-inc}) yields
\begin{equation}\label{fbeta} f^*(\beta) \cdot f^*\fitt_{k+1}(f_*\mathcal{O}_{X,0}) \subseteq (f^*\fitt_{k+2}(f_*\mathcal{O}_{X,0}))^2. \end{equation}
At this stage, we will refer to the assumption of the induction that (\ref{thm-fitk}) holds for $k$. For that we need a few results regarding $D^2(f)$ and $\pi^2_1$. 

Firstly, since $f$ is finitely $\mathcal{A}$-determined and of corank 1, double point space $D^2(f)$ is an isolated complete intersection singularity of dimension $n-1$ by \cite[Lemma 3.2]{conejero-nuno} (which is based on expanding \cite[Theorem 2.14(i)]{marar-mond} to the case of map-germs over ICIS). Moreover, $\pi^2_1: D^2(f) \rightarrow (X,0)$ is also finite and of corank 1 (\cite[Lemma 3.10]{klu-curv}); furthermore, $\pi^2_1$ is also finitely $\mathcal{A}$-determined since $D^{k-1}(\pi^2_1) = D^k(f)$ by the Principal of iteration (\cite[Lemma 3.9]{klu-curv}), and  $D^k(f)$ is an isolated complete intersection singularity of dimension $n+1-k$, or empty if $k>n+1$ (\cite[Lemma 3.2]{conejero-nuno}).
Therefore, the assumption of the induction for $k$ applied to $\pi^2_1$ yields
\begin{equation}\label{fitk-pi21}(\fitt_{k+1}^2((\pi^2_1)_*\mathcal{O}_{D^2(f),0}) : \fitt_{k}((\pi^2_1)_*\mathcal{O}_{D^2(f),0})) 
		= \fitt_{k+2}((\pi^2_1)_*\mathcal{O}_{D^2(f),0}).\end{equation}
Now, by the analytic version of the principal of iteration (see the proof of \cite[Lemma 3.9]{klu-curv}),
\[\fitt_{r}((\pi^2_1)_*\mathcal{O}_{D^2(f),0})=f^*\fitt_{r+1}(f_*\mathcal{O}_{X,0})\]
for any $r\geq 2$. Therefore,  (\ref{fitk-pi21}) yields
\begin{equation}\label{fitk-pi-f}
		((f^*\fitt_{k+2}(f_*\mathcal{O}_{X,0}))^2: f^*\fitt_{k+1}(f_*\mathcal{O}_{X,0})) 
		=f^* \fitt_{k+3}(f_*\mathcal{O}_{X,0}).
\end{equation}
So, from (\ref{fbeta}), we must have $f^*(\beta)\in f^* \fitt_{k+3}(f_*\mathcal{O}_{X,0})$, that is, $\beta \in\fitt_{k+3}(f_*\mathcal{O}_{X,0})$. This completes the proof.
\end{proof}

\begin{ex} Revisiting our example, we calculate the remaining Fitting ideals using Theorem \ref{thm2} and obtain the same results as in Example \ref{ex0}.

\begin{verbatim}
LIB "ring.lib";
ring t=0,(X,Y,Z),(wp(4,5,6));
ring s=0,(x,y),(wp(4,1));
ideal p=0;
map f=t,x,y5-xy,y6+xy2;
setring t;
ideal h=preimage(s,f,p);
ideal jh=jacob(h);
ideal fit1=radical(jh);
ideal fit2=std(quotient(fit1*fit1,h));
ideal fit3=std(quotient(fit2*fit2,fit1));
ideal fit4=std(quotient(fit3*fit3,fit2));

\end{verbatim}
\end{ex}

\begin{rem} Examples suggest that the statement of Theorem \ref{thm2} also holds for map-germs of corank 2. Interested readers may consider the map-germ
\begin{eqnarray} \label{ex-cor2} f\colon (\CC^3,0)&\rightarrow &(\CC^4,0)\\
(x,y)&\mapsto &(x,y^2+xz,z^2+xy,y^3+y^2z-yz^2+z^3)\nonumber \end{eqnarray}
which has $\mathcal{A}_e$-codimension 33 (\cite{altintas-ex2}) and $Q(f)=\mathbb{C}\cdot\{1,y,z,yz\}$ to calculate the Fitting ideals.  
\end{rem}


\appendix

\section{Calculating a presentation}\label{app1}
For comparison purposes, here we will calculate the Fitting ideals for our main example using Mond and Pellikaan's algorithm in \cite{mond-pellikaan} (also see  \cite[11.1]{mond-nuno}).

\begin{ex} \label{ex0} Let $f\colon (\CC^2,0)\rightarrow (\CC^3,0)$ be given by $f(x,y)=(x,y^5-xy,y^6+xy^2)$ as in (\ref{ex0-map}). We have $Q(f)=\CC\cdot\{1,y,y^2,y^3,y^4\}$. We find that
\begin{eqnarray*} Z=y^6+xy^2&=&Y\cdot y+2X\cdot y^2\\
Z\cdot y= y^7+xy^3&=&Y\cdot y^2+2X\cdot y^3\\
Z\cdot y^2= y^8+xy^4&=&Y\cdot y^3+2X\cdot y^4\\
Z\cdot y^3= y^9+xy^5&=&2XY+2X^2\cdot y+Y\cdot y^4\\
Z\cdot y^4= y^{10}+xy^6&=&Y^2+3XY\cdot y +2X^2\cdot y^2
\end{eqnarray*}
(Multiply the first relation by $-X$ and add to the last one to make it symmetric.) It follows that
\begin{equation}\label{ex0-lambda} \lambda = \left[\begin{array}{ccccc} 
-Z & 0 &   0&    2XY &Y^2 +XZ \\
Y &  -Z &  0 &   2X^2 &   2XY\\
2X & Y &    -Z & 0 &   0 \\
0 &   2X &  Y &    -Z & 0    \\
0 &  0 &    2X &   Y &    -Z 
\end{array} \right ]
\end{equation}
and that
\begin{eqnarray*}\fitt_0&=&(16X^5Y^2+Y^6-16X^6Z+11XY^4Z+28X^2Y^2Z^2+8X^3Z^3-Z^5)\\
\fitt_1&=& (16X^5-Y^4-6XY^2Z-4X^2Z^2,16X^4Y+Y^3Z+4XYZ^2, \\ &&
 8X^3Y^2+8X^4Z-Y^2Z^2-2XZ^3, 4X^2Y^3+12X^3YZ+YZ^3, \\ && 2XY^4+10X^2Y^2Z+4X^3Z^2-Z^4, Y^5+7XY^3Z+8X^2YZ^2) \\
\fitt_2&=& (X^3, X^2Y, X^2Z, XY^2, XYZ, Y^3, XZ^2, Y^2Z,YZ^2, Z^3) \\
\fitt_3&=& (X^2, XY, XZ, Y^2, YZ, Z^2) \\
\fitt_4&=&(X,Y,Z).
\end{eqnarray*}

The algorithm is implemented into a {\sc{Singular}} library \textsf{presmatrix.lib} by Hernandes, Miranda and Pe{\~n}afort-Sanchis in \cite{hernandes} and can be downloaded at \cite{presmatrix-lib}. Here we show a code for our example.

\begin{verbatim}
LIB "presmatrix.lib";
ring t=0,(X,Y,Z),(wp(4,5,6));
ring s=0,(x,y),(wp(4,1));
map f=t,x,y5-xy,y6+xy2;
presmatrix(f,0);
\end{verbatim}

Alternatively, a presentation matrix can be calculated in \textit{Macaulay2} (\cite{M2}) using the commands included in the program as follows.

\begin{verbatim}
S=QQ[x,y,Degrees=>{4,1}]
T=QQ[X,Y,Z,Degrees=>{4,5,6}]
f=map(S,T,{x,y^5-x*y,y^6+x*y^2})
pr=relations trim pushForward(f,S^1)
\end{verbatim}

\end{ex}

\section*{Acknowledgments}
The question of calculating Fitting Ideals without a presentation was first suggested to the first author by Prof. David Mond, to whom she is sincerely grateful for all he taught and passing the love of theory of differentiable map-germs. We send our warm gratitude to our families who supported us constantly over the course of this work.

\bibliographystyle{amsplain}
\bibliography{a-reference}

\end{document}